\documentclass[12pt]{amsart}

\usepackage[british]{babel}

\usepackage[latin1]{inputenc}

\usepackage{tabularx}
\usepackage{amscd}
\usepackage{times}
\usepackage[T1]{fontenc}

\newcommand{\PP}{\mathbb{P}}
\newcommand{\Gr}{\mathbb{G}}
\newcommand{\OO}{\mathcal{O}}

\newcommand{\HH}{\mathcal{H}}

\usepackage[titletoc,toc,title]{appendix}
\usepackage{amsfonts}
\usepackage{amsmath}
\usepackage{amsthm}
\usepackage{euscript}
\usepackage{amssymb}
\usepackage{enumerate}
\usepackage{hyperref}
\usepackage{mathrsfs}
\usepackage{xypic}
\usepackage{stmaryrd}

\xyoption{all}

\theoremstyle{plain}
\newtheorem{thm}{Theorem}[section]

\newtheorem{question}{Question}

\theoremstyle{definition}

\theoremstyle{remark}

\numberwithin{equation}{section}

\newcommand{\G}{\mathbb G}

\SelectTips {cm}{}

\begin{document}

\title[Geometry of lines and degeneracy loci]{Geometry of lines and degeneracy loci of morphisms of vector bundles}

\author
{Emilia Mezzetti}
\address{Dipartimento di Matematica e Geoscienze,
  Universit\`a degli studi di Trieste, Via Valerio 12/1,
34127 Trieste -- Italy\\
\url{mezzette@units.it}
}


\begin{abstract} Corrado Segre played a leading role in the foundation of line geometry. We survey some recent results on degeneracy loci of morphisms of vector bundles where he still is of profound inspiration.

\end{abstract}
\thanks{Acknowledgments:  The author is member of GNSAGA and is supported by FRA - Universit\`a di Trieste - and PRIN \lq\lq Geometria delle variet\`a algebriche''. }
\maketitle

\section{Grassmannians of lines: linear sections and focal properties}

\subsection{Classical point of view: work of Corrado Segre on Grassmannians}\label{classical}

The geometry of  families of lines in the projective space, classically called \lq \lq\thinspace line geometry'', has been one of the first objects of investigation of Corrado Segre, and a {\it fil rouge} of his research throughout his career. His graduation thesis was published in 1883, in two articles, in  \lq\lq\thinspace Memorie  dell'Accademia delle Scienze di Torino''. In the first article \cite{segre1883_1}, corresponding to the first two chapters of the thesis, Segre systematically studies the (hyper)quadrics,  in the second one \cite{segre1883_2} he develops the geometry of the Klein quadric in $\PP^5$, i.e. the Grassmannian $\Gr(1,3)$ of lines in $\PP^3$. 
He considers  linear and quadratic sections of the Grassmannian, called  linear and quadratic complexes when of dimension three, and  linear and quadratic congruences when of dimension two,  and moreover ruled surfaces. In particular he studies the notion of focal locus of a congruence, roughly speaking the set of points where two "infinitely near" lines of the family meet. 

We want to mention then a short article published by Corrado Segre in
1888 \cite{segre1888}, where he considers the line geometry in a projective space of any dimension $n$. Here Segre states a few basic facts on the
focal points of a  congruence of lines in $\PP^n$, meaning now a family of dimension $n-1$.

Many years later, in 1910, Corrado Segre comes back to this circle of ideas in the foundational article
\cite{segre1910}, entirely devoted to the 
projective-differential geometry of families of linear spaces, and in particular he develops the theory of their focal loci.

Moving in a different direction, in 1917 Segre publishes  the article \cite{segre1917}, where there is
a systematic study of the Grassmannian $\Gr(2,5)$ of $2$-planes in $\PP^5$, for instance
the classification of the orbits of the natural action of the projective group $PGL(6)$ on $\PP(\Lambda^3\mathbb C^6)$.
\medskip

An interesting article on the line geometry, possibly inspired by Corrado Segre, was published in 1891 by 
 Guido Castelnuovo \cite{castelnuovo}.

Here Castelnuovo studies the linear sections  of all dimensions, both general and special, of the Grassmannian of lines in $\PP^4$, $\Gr(1,4)$,  giving a geometric description of their  \lq\lq singular'' (i.e. focal) varieties in $\PP^4$. In the first cases they are
a conic, a projected Veronese surface (in the case of a congruence of lines), and
 a Segre cubic hypersurface.

This last hypersurface had been considered by Segre in the article \cite{segre_cubica}, which contains a description of its beautiful and rich geometry. The point of view developed by Castelnuovo allows to recover  the geometrical properties of this threefold with a completely different but natural method.  As Castelnuovo says, his article gives a
prototype of the geometry of lines in spaces of even dimension.

Let us look more closely at the particular case of a general linear congruence of lines  in $\PP^4$, let us call it $\Gamma$. In other words $\Gamma$ is a $3$-dimensional general linear section of $\Gr(1,4)$. Castelnuovo proves that
the elements of  $\Gamma$ are the trisecant lines of its focal variety, which results to be a projected Veronese surface $S$ in $\PP^4$. This implies that $\Gamma$ can be reconstructed from $S$. On the other hand, Castelnuovo proves that all projected Veronese surfaces in $\PP^4$ arise from such a construction.

To conclude this short historical overview, we want to mention that 
Francesco  Palatini in 1900 \cite{palatini} considered similar constructions in $\PP^5$. Here a general linear congruence gives raise, as singular variety, to an interesting three-dimensional smooth scroll $X$, after him named \lq\lq Palatini scroll''.
The interest comes from a \lq\lq non-normality'' property of this threefold, being in some sense the analogous of one of the characteristic properties of the Veronese surface. More precisely, the famous Severi theorem asserts that  the projected Veronese surface
 is the only non degenerate non linearly normal surface of $\PP^4$. 
It easily results from the construction that Palatini scroll  $X$ is smooth and not quadratically normal, meaning that the linear series cut on $X$ by the quadric hypersurfaces is not complete.  An important conjecture by Christian Peskine  states that this is the only smooth and non quadratically normal threefold in $\PP^5$ \cite{conj}.

\subsection{Modern  point of view}\label{modern}
In this section we will expose how the classical constructions mentioned in Section \ref{classical} can be rephrased in the modern language of the theory of vector bundles. The main references are \cite{ott-scrolls}, \cite{bazan-m},  \cite{dp}, \cite{tanturri1}. 

Let $k$ be an algebraically closed field, with char$(k)=0$. We fix integer numbers
$2\leq m\leq n$
and two $k$-vector spaces $U,V$, with $\dim U=m$, $\dim V=n+1$.

We consider a general bundle map 
$\phi: U^*\otimes \OO_{\PP(V)}\rightarrow \Omega_{\PP(V)}(2).$ 
From the Euler sequence, it follows that $H^0(\Omega_{\PP(V)}(2))\simeq (\Lambda^2V)^*$. Hence
$\phi$ is defined by $m$ general skew-symmetric matrices  $A_1,..., A_m$. Let
$X_\phi$ be the degeneracy locus of $\phi$. 
Then
\begin{itemize}
\item $\phi$ is  injective;
\item $\dim(X_\phi)=m-1$; 
\item $\dim(Sing(X_\phi))=2m-n-4;$
\item $X_\phi$ is set-theoretically the set of points $x\in \PP^n$ such that  $ (y_1A_1+\dots+y_mA_m)x=0$, for some $[y_1,...,y_m]\in\PP(U)$;  
\item for $n$ even, $X_\phi$ is a unirational variety parametrized by $\PP(U)$;
\item for $n$ odd, $X_\phi$ is a scroll over the hypersurface $Z$ in $\PP(U)$ defined by the equation: Pfaffian$(y_1A_1+\dots+y_mA_m)=0$.

\end{itemize}
The skew-symmetric matrices  $A_1,..., A_m$ correspond bijectively to  hyperplanes $H_1, ... , H_n$ in the embedding space of the Grassmannian $\Gr(1,n)$, so they define a linear section $\Gamma:=\Gr(1,n)\cap H_1\cap...\cap H_n$, which is nothing else than  the family of lines of the classical construction, and the degeneracy locus $X_\phi$ coincides with the classical focal locus.

For $n=4$, this construction is precisely the one described by Castelnuovo. If $m=2$ $X_\phi$ is a plane conic, if $m=3$ a projected Veronese surface, if $m=4$ a cubic Segre hypersurface.
For $n=5$, if $m=3$ we get an elliptic surface scroll  in $\PP^5$, if $m=4$ a Palatini scroll. 

Having in mind the case of the Veronese surface in $\PP^4$, the following natural question arises.
\begin{question}  Can $\Gamma$  be recovered from $X_\phi$?\label{q1}
\end{question}
 The example of an elliptic surface scroll  in $\PP^5$ shows that the answer is in general negative. Indeed such a surface is obtained as degeneracy locus from four different maps $\phi$. This had been observed by Gino Fano in \cite{fano}  (see \cite{bazan-m} for a modern proof), and follows from the structure of the Picard group of elliptic curves.  

On the other hand,  a Palatini scroll $X$
 is obtained from only one map $\phi$
\cite{fania-m1}. \ The geometrical interpretation of the congruence $\Gamma$ in this case is that the lines of $\Gamma$ are the $4$-secant lines of $X$ not contained in $X$, so unicity follows.


 \subsection {General linear sections: Hilbert schemes}
In this section we reformulate Question \ref{q1} in more precise terms, and we describe the recent progress on it by D. Faenzi- M.L. Fania \cite{faenzi-fania} and F. Tanturri  \cite{tesiT}, \cite{tanturri2}, \cite{tanturri1}.

Let us denote by
  $\mathcal H$  the union of the irreducible components of the Hilbert scheme containing the degeneracy loci of general maps 
$\phi: U^*\otimes \OO_{\PP(V)}\rightarrow \Omega_{\PP(V)}(2).$

We consider the natural rational map $\rho: \Gr(m, \Lambda^2V)\dashrightarrow \HH$, taking a bundle map to its degeneracy locus.
The discussion above motivates the following:

\begin{question}\label{q2} Let $n$, $m$ be integer numbers with $2\leq m\leq n+1$. 
\begin{enumerate}
\item Is $\rho$ dominant? 
\item Is $\rho$ generically injective? 
\item If $\rho$ is not injective, describe its fibres.
\end{enumerate}
\end{question}

The results of the above quoted authors give a complete answer to Question \ref{q2} for $2\leq m\leq n$. They are the content of the following theorem.


\begin{thm}  In the above notation,
\begin{itemize}
\item if $m\geq 4$ or $(m,n)=(3,4)$, then $\rho$ is birational and $\HH$ is generically smooth;
\item if $(m,n)=(3,5)$, then $\rho$ is dominant and $4\colon1$;
\item if $m=3$ and $n\neq 5$, then $\rho$ is generically injective;
\item if $m=3$ and $n$ is even,  a general element of the image of $\rho$ is a special projection of the Veronese surface $v_{\frac{n-1}{2}}(\PP^2)$. The centre of projection is  the linear span of the partial derivatives of order $\frac{n-5}{2}$ of a non-degenerate polynomial of degree $n-3$;
\item if $m=3$ and $n$ is odd, a general element of the image is a projective bundle $\PP(\mathcal G)$, with $\mathcal G$ a general stable rank two bundle on a general plane curve $C$ of degree $\frac{n}{2}$,  with $\det(\mathcal G)=\OO_C(\frac{n-2}{2})$;
\item if $m=2$, then $\rho$ is dominant. Its fibres have positive dimension and can be explicitly described.
\end{itemize}
\end{thm}

The proof of the theorem relies on 
the method of Kempf - Lascoux - Weyman, in particular on applications of Eagon--Northcott complex, and on the structure theorem of Buchsbaum-Eisenbud and
apolarity theory  in the case $m=3$.

\section{Linear systems of skew--symmetric matrices of constant rank}

Up to now we considered  {\it general}  linear sections of Grassmannians of lines, or in other words, general bundle maps $\phi: U^*\otimes \OO_{\PP(V)}\rightarrow \Omega_{\PP(V)}(2).$ But, according to a general idea by Ch. Peskine and F. Zak, we expect that  degeneracy loci of {\it special} maps $\phi$ should produce interesting varieties.
In this section we will explore some special cases and related problems.

As we observed in Section \ref{modern}, giving the map $\phi$ 
is equivalent to giving an $m$-dimensional  linear subspace $\Lambda:=\langle A_1,..., A_m\rangle$ in $\Lambda^2 V^*$.
The family of lines $\Gamma$ is the intersection of $\Gr(1,n)$ with $\PP(\Lambda^*)$, the projectivized dual of $\Lambda$.
Hence, 
classifying the maps $\phi$, for fixed $m$, is equivalent to classifying  the projective subspaces of dimension $m-1$ in $\PP(\Lambda^2 V)^*$ under the action of $PGL(n+1)$, also describing their possible positions with respect to the natural filtration by the rank of skew-symmetric tensors. 

Since a complete classification appears to be  out of reach, we will focus on a natural special case, which is interesting  and has been studied  in linear algebra since the classical work of Gantmacher \cite{Gantmacher}.

By $\sigma_k\Gr(1,n)$ we denote the $k$th secant variety, i.e. the Zariski closure of the union of the linear spans of general $k$-tuples of points in $\Gr(1,n)$.

\begin{question}\label{const rank}
Classify  the orbits of  linear subspaces contained in a quasi-projective variety of the form $\sigma_k\Gr(1,n)\setminus \sigma_{k-1}\Gr(1,n)$ for some $n$ and $k$, i.e. the orbits of linear systems of skew-symmetric matrices of constant  rank.
\end{question}

This question is of course a particular case of the more general problem of classifying linear systems of matrices of constant rank of any size without skew--symmetry conditions. An analogous problem refers to symmetric matrices. A different generalization regards linear systems of matrices of bounded rank. For   a  bibliography on these problems and their applications in algebraic geometry, see \cite{Ilic_JM}.

From now on, when speaking of dimension of a linear system we will mean its {\it projective} dimension.

The linear systems of skew-symmetric matrices of constant rank $2$ can be interpreted as linear spaces contained in a Grassmannian $\G(1,n)$, a  well-known case. Those of maximal dimension correspond either to the lines contained in a $2$-plane or to a star of lines; they are a $\PP^2$ and a $\PP^{n-1}$ respectively. So the first interesting  case to analyze is that of linear systems of $6\times 6$ skew-symmetric matrices of constant rank $4$. This situation has been studied in \cite{Manivel_Mezzetti}. The result is the following classification: 
\begin{itemize}
\item the maximal dimension of  a linear system of skew--symmetric $6\times 6$ matrices of constant rank four is two, 
\item the space of $1$-dimensional linear systems of matrices of this type is irreducible of
dimension $22$, with an open $PGL(6)$-orbit of general lines, and a codimension
one orbit of special lines,
\item  there are four $PGL(6)$-orbits of $2$-dimensional linear systems, all of the same dimension $26$ and homogeneous under the action of $PGL(6)$.
\end{itemize}

The above classification of the orbits of $2$-planes  has found application to the classification of the 
degenerations of  Palatini scroll \cite{DP_M},  such that the cubic surface, which is the base of the scroll, splits as union of a $2$-plane and a quadric surface.  As a further application, we mention the construction of new examples of non quadratically normal threefolds in $\PP^5$ \cite{dpm_quadr}.  These threefolds are singular, so the smoothness assumption in the conjecture of Peskine mentioned in Section \ref{classical} results to be necessary.


\subsection{Maximal dimension}

For linear systems of $(n+1)\times (n+1)$ skew-symmetric matrices of constant rank $r$, with  $r<n+1$ even, the first step in the study of Question \ref{const rank} is finding  their maximal dimension.

Let us introduce the following notation:
$l(r,n+1):=\max \{\dim$ of a linear system of  skew-symmetric matrices of rank $r$ and size $n+1\}.$

The following inequality is due to J. Sylvester \cite{sylvester} and R. Westwick \cite{Westwick1}:  if $r\geq 2$ is an even number, 
$r<n+1$, then 
\begin{equation}\label{disug}
n+1-r\leq l(r,n+1) \leq 2(n+1-r).\end{equation}
In particular, if $r=n-1$, i.e. if the corank is two, then $2\leq l(n-1,n+1) \leq 4$. 

We mention that the same inequality is valid also for linear systems of symmetric matrices.
The main result in \cite{Ilic_JM} states that  the maximal dimension of a linear system of symmetric matrices of constant rank is $n+1-r$, the lower bound above.
We will see that in the skew-symmetric case the situation is different, in particular $l(r,n+1)$ depends on $n$, and not only on $n-r$, and its precise value is still unknown in general.

To explain why and to underline the deep connections of this problem with algebraic geometry, let us remark that  a linear system of  dimension $d\leq l(r,n+1)$ of skew-symmetric matrices of constant rank defines an exact sequence of vector bundles:


\begin{equation}\label{ex_seq} 0 \rightarrow K\rightarrow \OO_{\PP^{d}}^{n+1} \rightarrow \OO_{\PP^{d}}^{n+1}(1) \rightarrow N\rightarrow 0,
\end{equation}
where:
\begin{itemize}
\item the kernel $K$ and the cokernel $N$ are rank $r$ vector bundles, 
\item $N\simeq K^*(1)$,
\item $K^*$ is globally generated and defines an embedding in the appropriate Grassmannian.
\end{itemize}

In particular, if $r=n-1$, then $K$ and $N$ are rank $2$ bundles and
$c_1(K^*)=\frac{r}{2}$.
For instance, for $n=5$ and $r=4$, the four orbits of $2$-planes of $6\times 6$ matrices of constant rank $4$ correspond to the four rank $2$ bundles  on $\PP^2$ with $c_1=2$ defining an embedding in $\Gr(1,5)$, according to their classification given in \cite{sierra_ugaglia_grassm}. 

This can be restated introducing the definition of  $m$-effective bundle. We say that a globally generated rank $r$ bundle on $\PP^d$ is $m$--effective if its dual appears as kernel $K$ in an exact sequence of the form (\ref{ex_seq}). So all globally generated rank $2$ bundles on $\PP^2$, with $c_1=2$ and defining an embedding in $\Gr(1,5)$, are $m$-effective.

A similar result is true for linear systems of
$8\times 8$ skew--symmetric matrices of constant rank $6$. In \cite{Fania_Mezzetti} it is proved that every globally generated rank $2$ bundle on $\PP^2$, with $c_1=3$ and defining an embedding in $\Gr(1,7)$, is $m$-effective. This gives a picture of the $2$--planes.  In the same article it is proved that also for these size and rank there are no linear systems of  dimension bigger than $2$. 

\subsection{Spaces of dimension $1$ and $2$}

Turning to general $n$, we  mention that in \cite{Fania_Mezzetti} a complete classification  of the orbits of $1$--dimensional spaces is given,  for every pair  $(n, r)$.

As for $2$-planes, we focus on the case $n$ odd and $r=n-1$, corresponding to  rank $2$ bundles on $\PP^2$. The following result has been proved in \cite{bo_me}.
\begin{thm} \label{2planes}
Let $n$ be an odd integer number.
\begin{itemize}
\item There exists an $m$-effective bundle $E$ with Chern classes $c_1, c_2$ for every pair $(c_1, c_2)$ such that $c_1=\frac{n-1}{2}$ and $c_2$ is in the stable range; 
\item there exist unstable non $m$-effective globally generated rank 2 bundles.
\end{itemize}
\end{thm}

The proof of the theorem rests on the recent complete description of the possible Chern classes of rank two globally generated vector bundles on the projective plane by Ph. Ellia \cite{ellia_effective}.

The assumption that $r=n-1$, or that the corank is two, is not so restrictive, because every $2$--plane of skew--symmetric matrices of constant rank $r$ can be isomorphically projected to a $2$--plane of $(r+2)\times (r+2)$ matrices (for more details see \cite[Cor. 5.9]{Fania_Mezzetti}). 


\subsection{3-planes of skew-symmetric matrices of corank two} 
According to Theorem \ref{2planes}, there are \lq\lq\thinspace many'' orbits of $2$-planes of skew-symmetric matrices of constant corank two. But the upper bound on $l(n-1, n+1)$ given by  inequality (\ref{disug}) is $4$. The first example of a linear system of dimension bigger than two was exhibited by Westwick in \cite{westwick}.  It is 
a $3$--space of $10\times 10$ skew-symmetric matrices of constant rank $8$.

Since this example is given without any explanation, it is natural to ask  which are the corresponding bundles and to look for  other examples. 

These problems have been considered in \cite{bo_fa_me}.
The main result is a method to contruct a skew-symmetric matrix $A$ of linear forms on $\PP^3$,
having constant rank $r=n-1$, i.e. corank $2$,
 starting from a normalized rank $2$ bundle $E$ on $\PP^3$ with  allowed Chern classes, i.e. $c_1(E)=0$ and $c_2(E) = \frac{r(r+4)}{48}$. 
The matrix $A$ will  have $E$ as kernel, up to a twist by a line bundle.
One  needs also a class
$\beta \in Ext^2(E(\frac r4 -1),E(-\frac r4 -2))$ .

Then the  main idea is the following: 
  the cone of $\beta$, interpreted as
a morphism $E(\frac r4 -1) \to E(-\frac r4 -2)[2]$
in the derived category $D^b(\PP^3)$, is a $2$-term complex. Using Beilinson's theorem,  necessary and sufficient conditions are found ensuring that this $2$-term complex is of the desired form (\ref{ex_seq}) and the matrix in the middle is skew--symmetrizable.

A useful remark is that 
the conditions can be 
simplified if $E$ has natural cohomology.
Since general instantons have natural cohomology,
and their minimal graded free
resolution is known, as well as that of their cohomology module,  it is possible to apply the construction to general instantons. It turns out that the necessary and sufficient conditions are satisfied for low values of $c_2(E)$.

The result is the following:
\begin{thm}

\begin{enumerate}
\item Any $2$--instanton on $\PP^3$ induces a $3$--dimensional space of $10 \times 10$ skew--symmetric matrices of
 constant  rank $8$;
\item any $4$-instanton $E$ on $\PP^3$ with natural cohomology and such that $E(2)$ is globally generated induces a $3$--dimensional space of $14 \times 14$ skew--symmetric matrices of
  constant rank $12$;
\item Westwick's example corresponds to a $2$--instanton belonging to the most special orbit of the moduli space $M_{\PP^3}(2;0,2)$ under the natural action of $SL(4)$.
\end{enumerate}
\end{thm}

Hence in particular
  there exists a continuous family of examples of $3$-spaces of $10 \times 10$ skew--symmetric matrices of
  rank $8$, all non--equivalent to Westwick's one.

Explicit constructions of all these matrices of linear forms can be found in \cite{bo_fa_le}.

\subsection{Open problems} 
To end this survey we list a few of the many  open problems related to linear systems of matrices of constant rank.

\begin{itemize}

\item Do there exist other examples of $3$-spaces of skew--symmetric matrices of
  constant corank $2$, besides those given by $2$-instantons and $4$-instantons? 

\item Do there exist linear systems of dimension $4$ of skew--symmetric matrices of
  constant corank $2$?

\item What kind of variety is the union of the orbits of spaces of matrices of fixed constant rank? 
\end{itemize}

\bibliographystyle{amsalpha}
\bibliography{Mezzetti_Segre}

\end{document}